\documentclass[12pt]{amsart}
\usepackage[mathscr]{eucal}
\usepackage[all]{xy}

\newtheorem{theorem}{Theorem}[section]
\newtheorem{lemma}[theorem]{Lemma}
\newtheorem{proposition}[theorem]{Proposition}
\newtheorem{corollary}[theorem]{Corollary}
\newtheorem{defn}[theorem]{Definition}

\theoremstyle{remark}
\newtheorem{remark}[theorem]{Remark}

\numberwithin{equation}{section}

\title{Quasi-Hodge Metrics and Canonical Singularities}
\author{Chin-Lung Wang}
\address{National Tsing-Hua University, Hsinchu, Taiwan}
\email{dragon@math.nthu.edu.tw}

\begin{document}
\maketitle

\section{Introduction and Statement}

Throughout this paper we work over the base field $\mathbb{C}$. The basic problem in
algebraic geometry treated here is the {\it filling-in problem} (or {\it
degeneration problem}). Given a smooth projective family
$\mathscr{X}^\times\to\Delta^\times$, we would like to know when we can fill in a
reasonably nice special fiber $\mathscr{X}_0$ to form a projective family
$\mathscr{X}\to \Delta$, perhaps up to a finite base change on $\Delta^\times$. For
example, when can $\mathscr{X}_0$ be smooth? When can it be irreducible? Or when can
it be irreducible with at most certain type of mild singularities? We would like to
search for conditions depending only on the punctured family.
\par

In this paper, for any smooth projective family $\mathscr{X}\to S$ over a
smooth base $S$ such that $\mathscr{X}_s$ has semi-ample canonical bundle, we
shall define for each large $m\in \mathbb{N}$ a K\"ahler metric $g_m$ on $S$,
called the $m$-th quasi-Hodge metric.  When $S=\Delta^\times$, we propose that
the incompleteness of $g_m$ near $0$ for suitable $m$'s provides a necessary
and sufficient condition for the existence of $\mathscr{X}_0$ to be irreducible
and with at most {\it canonical singularities} (c.f.\ Remark 2.5). Notice that
the metric incompleteness condition is insensitive to base changes.\par

More precisely, for a smooth projective family $\pi:\mathscr{X}\to S$ with
$p_g(\mathscr{X}_s)\ne 0$, the (possibly degenerate) quasi-Hodge metric
$g_H=g_1$ on $S$ is given by the semi-positive first Chern form of the rank
$p_g$ Hodge bundle $F^n = \pi_*K_{\mathscr{X}/S}$. When $S = \Delta^\times$,
let $T$ be the monodromy operator acting on $H^n(\mathscr{X}_s,\mathbb{C})$
where $n=\dim\mathscr{X}_s$, $s\ne 0$.

\begin{theorem}
For any smooth projective family $\pi:\mathscr{X}\to \Delta^{\times}$,
\begin{itemize}
\item[(1)] The quasi-Hodge metric $g_H$ is incomplete if and only if that
$H^{n,0}(\mathscr{X}_s)$ consists of $T^r$-invariant cycles for some
$r\in\mathbb{N}$.
\item[(2)] In terms of a semi-stable model with
$\mathscr{X}_0=\bigcup_{i=0}^N X_i$, $g_H$ is incomplete if and only if that
$p_g(\mathscr{X}_s)=\sum_{i=0}^N p_g(X_i)$.
\item[(3)] In particular, degenerations with Gorenstein canonical singularities
have finite $g_H$ distance.
\end{itemize}
\end{theorem}

This generalizes an earlier result on Calabi-Yau families in \cite{wang} concerning
Weil-Petersson metrics (c.f.\ \S7). In the Calabi-Yau case, since $p_g=1$, there is
exactly one component $X_i$ with $p_g\ne 0$. Starting from this we may deduce

\begin{proposition}
Let $\mathscr{X}\to \Delta$ be a degeneration of Calabi-Yau $n$-folds such that
$0\in\Delta$ has finite $g_H$ distance. Then the minimal model conjecture (MMC)
in dimension $n+1$ implies that, up to a finite base change, $\mathscr{X}\to
\Delta$ is birational to $\mathscr{X}'\to \Delta$ such that
$\mathscr{X}_s\cong\mathscr{X}'_s$ for $s\ne 0$ and $\mathscr{X}'_0$ is a
Calabi-Yau variety with at most canonical singularities.
\end{proposition}

This had been solved in the literature for abelian varieties and K3 surfaces
(type I degenerations in Kulikov's work \cite{kulikov}). Indeed in both cases
one ends up with smooth families. Here we show that the relative minimal model
is the filling-in we seek for. The primary reason is that all other components
in $\mathscr{X}_0$ are uni-ruled hence should be contractible. The MMC is known
in dimensions up to three by the Mori theory. Recently Shokurov announced the
existence of log-flips in dimension four. Its validity would imply the result
for degenerations of Calabi-Yau threefolds.
\par

For the general cases when $p_g(\mathscr{X}_s)> 1$, there could be more than
one {\it essential components} with $p_g(X_i)\ne 0$, so $g_H$ is insufficient
for our purpose. In order to concentrate all the pluri-genera $P_m$ in one
component, we need a finer metric on the punctured disk. The actual
construction in \S3 is to consider the case that $\mathscr{X}_s$ has {\it
semi-ample canonical bundle} and to take an $m$-cyclic cover $\mathscr{Y}\to S$
of $\mathscr{X}\to S$ along a smooth divisor $\mathscr{D}\in
|K_{\mathscr{X}/S}^{m}|$ for suitable $m$'s. We then define the $m$-th
quasi-Hodge metrics $g_m$ to be $g_H$ with respect to $\mathscr{Y}\to S$. The
semi-ample assumption is to guarantee the invariance of pluri-genera in
projective families, which is not needed in Theorem 1.1 by the well-known
invariance of Hodge numbers.

\begin{proposition}
Let $\pi:\mathscr{X}\to\Delta$ be a degeneration of smooth projective
manifolds $\mathscr{X}_s$ with semi-ample canonical bundle. If
$\mathscr{X}_0$ is irreducible with only
canonical singularities then $0\in \Delta$ is at finite
$g_m$ distance for all $m$ such that $g_m$ is defined.
\end{proposition}

We conjecture that the converse statement holds true. To support some evidence,
we verify it in \S4 for degenerations of curves.

\begin{theorem}
Let $\pi:\mathscr{X}^\times\to \Delta^\times$ be a projective family of smooth
curves of genus $\ge 2$, then $g_m$ is defined for all $m\in\mathbb{N}$ and the
incompleteness of $g_m$ for any three values of $m$'s implies that up to a
finite base change $\pi$ can be completed into a smooth family.
\end{theorem}

In \S6 we discuss the notion of {\it essential incomplete boundary points}.
Namely a finite distance degeneration without any smooth filling-in up to any
finite base change. Notice that there is no such essential degenerations for
curves (Theorem 1.4), abelian varieties and K3 surfaces (as mentioned
above).\par

\begin{theorem}
There exist essential finite $g_m$ distance degenerations for threefolds:
\begin{itemize}
\item[\rm (1)] Nodal degenerations of Calabi-Yau threefolds with
$h^1(\mathscr{O})=0$: these provide the simplest type of examples with nontrivial
monodromy.
\item[\rm (2)] Terminal degenerations of smooth minimal threefolds: some of
them have trivial $C^\infty$ monodromy hence provide more subtle examples.
\end{itemize}
\end{theorem}

In \S7 we make a few remarks on the Weil-Petersson metric on Calabi-Yau moduli
spaces and its relation to Viehweg's theory on moduli spaces.\par

\subsection*{Acknowledgement.}
Part of this work was done during my visiting of Harvard University in the
spring of 2001. I am grateful to Professor S.-T.\ Yau and Professor W.\ Schmid
for their encouragement. I am also grateful to Professor M.\ Gross for useful
discussions during his visiting of NCTS, Taiwan in 2001.

\section{Quasi-Hodge Metrics $g_H$}

Let $\pi:\mathscr{X}\to S$ be a family of polarized algebraic manifolds
with $H^{n,0}(\mathscr{X}_s)=H^0(\mathscr{X}_s,K_{\mathscr{X}_s})\ne 0$.
It is well known that these spaces have constant rank $p_g$ (the geometric
genus) in $s\in S$ and they form a holomorphic vector bundle $F^n =
\pi_*K_{\mathscr{X}/S}$. It is the last piece of the Hodge filtration
$F: F^0\supset F^1\supset\cdots\supset F^n$
with a natural Hermitian metric induced from the topological cup product on each fiber:
$$
Q(u,v)|_{F_s}=\sqrt{-1}^n\int_{\mathscr{X}_s}u\cup v.
$$

Griffiths has shown that the first Chern form of $(F^n, Q)$ (for a local frame
$\{u_i\}$): $$ \omega=\frac{\sqrt{-1}}{2}\,{\rm Ric}_{Q}(H^{n,0})
    =-\frac{\sqrt{-1}}{2}\,\partial\bar\partial\log\,\det
    \Big[Q(u_i,\bar{u}_j)\Big]_{i,j=1}^{p_g}
$$
is semi-positive, hence defines a (possibly degenerate) K\"ahler metric on $S$.
We call it the quasi-Hodge metric. For Calabi-Yau families this agrees with
the Weil-Petersson metric defined via the Ricci flat metric on each fiber, but
not in other cases (c.f.\ \S7).\par

When $\pi:\mathscr{X}\to\Delta$ is a degeneration, i.e.\ $\pi$ is smooth outside the
puncture, we are interested in relating the differential-geometric properties of $\Delta$
under the quasi-Hodge metric and the singularities occur in $\mathscr{X}_0$.
Applying Mumford's semi-stable reduction theorem we may and will first assume that
$\mathscr{X}$ is smooth and $\mathscr{X}_0$ is a simple normal crossing divisor
in it. In this case the monodromy $T$ is unipotent. Let $N=\log\,T$ be the
nilpotent monodromy associated to it. We need Schmid's theory of
limiting mixed Hodge structure \cite{schmid} to analyze this situation. \par

Fix a reference fiber $X=\mathscr{X}_s$ with $s\ne 0$ and let $V$ to be
the primitive cohomology of $H^n(X,\mathbb{C})$.
Recall that $\pi$ induces a variation of Hodge structures (VHS) of weight $n$
on $V$ over $\Delta^{\times}$ and gives
rise to the period map $\phi:\Delta^{\times} \to \langle T\rangle\backslash D$,
with $D$ the period domain. The map $\phi$ lifts to the upper half plane
$\Phi:\mathbb{H}\to D$ with the coordinates $s\in\Delta^\times$ and $z\in\mathbb{H}$
related by $s=e^{2\pi\sqrt{-1} z}$. Set
$$
A(z)=e^{-zN}\Phi(z):\mathbb{H}\to\check{D}.
$$
($\check{D}$ is the compact dual of $D$.) Since $A(z+1)=A(z)$,
$A$ descends to a function $\alpha(s)$ on $\Delta^\times$. Schmid's
Nilpotent Orbit Theorem implies that $\alpha(s)$ extends holomorphically
over $s=0$. The special value $F_\infty:=\alpha(0)$ is
called the limiting filtration. \par

$N$ uniquely determines the monodromy weight filtration on $V$:
$0\subset W_0\subset W_1\subset\cdots\subset W_{2n-1}\subset W_{2n}=V$
such that $N(W_k)\subset W_{k-2}$ and induces an isomorphism
$$
N^\ell:G^W_{n+\ell}\cong G^W_{n-\ell},
$$
on graded pieces. $F^p_\infty$ and $W_k$ together define a
polarized mixed Hodge structure (MHS) on $V$: namely
the induced Hodge filtration
$$
F^p_\infty G^W_k:=F^p_\infty\cap W_k/F^p_\infty\cap W_{k-1},
\quad p=0,\ldots, n
$$
is a pure Hodge structure of weight $k$ on $G^W_k$. $N$ acts on them as a morphism
of type $(-1,-1)$ --- $N(F^p_\infty G^W_k)\subset F^{p-1}_\infty G^W_{k-2}$.
Moreover, for $\ell\ge 0$, the primitive part $P^W_{n+\ell}:=\ker
N^{\ell+1}\subset G^W_{n+\ell}$ is polarized by $Q(\cdot,N^\ell\bar\cdot)$.
By adding the non-primitive part, the total cohomology
$H^n(X,\mathbb{C})$ admits non-polarized MHS.\par

In the rest of this section we prove Theorem 1.1 by dividing it into
steps. We start with the following theorem which generalizes the one in \cite{wang}:

\begin{theorem}
The center $0\in\Delta$ is at finite geodesic distance from the generic point
$s\ne 0$ under the quasi-Hodge metric if and only if $NF^n_\infty=0$.
\end{theorem}

\begin{proof}
Let $\Phi:\mathbb{H}\to D$ be the lifting of the period map.
To start the computation, we need to choose a good holomorphic frame $\Omega_j$,
$j=1,\ldots,p_g$ of $F^n$. Let
$p^n: D\to G(p_g,V)$ be the projection to the $F^n$ part. we have
$\Phi^n(z)=(e^{zN}\alpha(s))^n=e^{zN}\alpha^n(s)$. Here
$*^n:=p^n(*)\in G(p_g,V)$ is the $n$-th flag.
Near $t=0$, we can represent $\alpha^n$ through
local homogeneous coordinates as $p_g$ vectors ${\bf a}^j$, $j=1,\ldots,p_g$
in $V$. Then ${\bf a}^j(s)=a^j_0 + a^j_1 s + \cdots$ is holomorphic in
$s$. We have correspondingly
$$
{\bf A}_j(z)=a^j_0 + a^j_1 e^{2\pi\sqrt{-1} z} + a^j_2 e^{4\pi\sqrt{-1} z} + \cdots.
$$

As in \cite{wang}, the function
$e^{2\pi\sqrt{-1} z}=e^{2\pi\sqrt{-1} x}e^{-2\pi y}$ has the
property that all the partial
derivatives in $x$ and $y$ decay to $0$ exponentially as $y\to\infty$, with
rate of decay independent of $x$. Let $h$ be the class of functions with this
property and ${\bf h}$ the corresponding class of functions with value in $V$.\par

Let $ \Omega_j(z)=e^{zN}{\bf A}_j(z) $ for $j=1,\ldots,p_g$.
This is the desired frame because frame representations correspond to
sections of the universal rank $p_g$ subbudle of $G(p_g,V)$ which pulls
back to $F^n$ by $\Phi$. So the K\"ahler form $\omega$ of the
induced quasi-Hodge metric on $\mathbb{H}$ is given by
$$
\omega=-\frac{\sqrt{-1}}{2}\partial\bar\partial\log\,
\det\Big[Q\big(e^{zN}{\bf A}_i(z),e^{\bar{z}N}\overline{{\bf A}_j(z)}\big)
\Big]_{i,j=1}^n.
$$

Since the base is one dimensional, if we write the metric as $G|dz|^2$ then
$ G=-(1/4)\triangle \log\,\det\,Q$. From
$Q(Tu,Tv)=Q(u,v)$, it follows easily that $Q(Nu,v)=-Q(u,Nv)$ and
$Q(e^{zN}u,v)=Q(u,e^{-zN}v)$. Since ${\bf A}_j=a^j_0+{\bf h}$, we have
\begin{eqnarray*}
\det\big[Q(e^{zN}{\bf A}_i,e^{\bar zN}\bar{\bf A}_j)\big]
\!\!\!&=&\!\!\! \det\big[Q(e^{zN}a^i_0,e^{\bar zN}\bar a^j_0) + h\big] \\
      &=&\!\!\! \det\big[Q(e^{2\sqrt{-1}yN}a^i_0,\bar a^j_0)\big] + h \\
      &=&\!\!\! p(y) + h,
\end{eqnarray*}
where $p(y)$ is a polynomial in $y$. Let $d=\deg\,p(y)$. It has the property that
$d=0$ if and only if $NF^n_\infty = 0$. This is a consequence of the polarization
condition for the mixed Hodge structure.\par

To see this, we may choose the basis ${\bf a}^j$ in such a way that
for $a^j_0\in G_{n + \ell_j}$ in the limiting MHS, $i<j$ implies
$\ell_i\ge\ell_j$. Let $\ell=\ell_1=\cdots=\ell_q>\ell_{q+1}$, that is
$\ell = \max\,\{\,i\in\mathbb{N}\cup\{0\}\,|\,N^i F^n_\infty\ne 0\,\}$.
Then the determinant of the $q\times q$ matrix corresponds to those entries with
index between $1$ and $q$ is a positive polynomial in $y$ of degree $q\ell$
by the polarization condition on $G_{n+\ell}$. Inductively we may find numbers
$q_k$ and $\ell_{(k)}$ with $\sum_k q_k = p_g$ and $\ell_{(k)}$ decreases such that
the determinant of the corresponding $q_k\times q_k$ matrix
(as a block on the diagonal) is a positive polynomial in $y$ of degree
$q_k\ell_{(k)}$. Now we conclude that the original determinant is a positive
polynomial in $y$ of degree $d=\sum_k \ell_{(k)}$ dominant by the product of
these diagonal blocks --- because all other elements have smaller degree in
$y$. It is clear that $d=0$ if and only if $\ell=\ell_1=0$ if and only if
$NF^n_\infty = 0$.\par

Then
\begin{eqnarray*}
4G \!\!\!&=&\!\!\! \frac{(p'+h)^2-(p+h)(p''+h)}{(p+h)^2}
              =\frac{(p'^2-p p'')+h}{p^2+h}\\
& & \!\!\!\sim \frac{p'^2-p p''}{p^2} + h
              \sim \frac{d^2-d(d-1)}{y^2} + h
              =\frac{d}{y^2} + h.
\end{eqnarray*}
Here we have used the fact that $p^{-2}h\in h$. Obviously, if
$NF^n_\infty=0$ then $d=0$ and $G=h$, so $\int^\infty_{z_0}\sqrt{G}\,|d z|<\infty $
for some curve (e.g. $x=c$).
When $NF^n_\infty\ne 0$ we have $d\ge 1$ and for $y$ large enough we
can make $h<1/y^3$
uniformly in $x$, then clearly
$ \int^\infty_{z_0}\sqrt{G}\,|d z|\sim 2\log y\,|^\infty_{y_0}=\infty $
for any path with $y\to\infty$.
\end{proof}

\begin{remark}
From the proof, we know that in the case of infinite distance,
the quasi-Hodge metric is exponentially asymptotic to a scaling of the
Poincar\'e metric.
\end{remark}

So far it is purely Hodge-theoretic and makes perfect sense in the framework of
abstract variation of Hodge structures. Now we plug in the geometric data:

\begin{theorem}
For a semi-stable degeneration $\pi:\mathscr{X}\to \Delta$ of varieties
with $p_g > 0$, if $\mathscr{X}_0 =\bigcup_{i=0}^N X_i$ then the induced
quasi-Hodge metric is incomplete at $0$ if and only if for $s\ne 0$,
$p_g(\mathscr{X}_s)=\sum_{i=0}^N p_g(X_i)$.
\end{theorem}

\begin{proof}
The proof is essentially the same as in \cite{wang}, so we only sketch it
briefly. \par

Deligne has shown that the cohomologies of the normal crossing divisor
$\mathscr{X}_0$ admit mixed Hodge structures. Let
$i:\mathscr{X}_s\to\mathscr{X}$ be the inclusion map and $i^\#:
H^n(\mathscr{X}_0)\cong H^n(\mathscr{X})\to H^n(\mathscr{X}_s)$ the induced
map. The Clemens-Schmid exact sequence
$$
\cdots\to H^n(\mathscr{X}_0)\mathop{\to}^{i^\#}
   H^n(\mathscr{X}_s)\mathop{\to}^N H^n(\mathscr{X}_s)\to\cdots
$$
is an exact sequence of MHS's. The is also known as the Invariant Cycle Theorem
which implies that $p_g:=p_g(\mathscr{X}_s)\ge\sum_{i=1}^N p_g(X_i)$ with
$\sum_{i=1}^N p_g(X_i)$ corresponds to $T$-invariant cycles in $F^n_s$. Now the
condition $NF^n_\infty=0$ says that all those $p_g$ cycles are $T$-invariant.
\end{proof}

Notice that in the statement of Theorem 1.1, part (1), $r$ is the degree of the
base change $t\mapsto t^r$ performed on $\Delta$ in obtaining the semi-stable
model.

\begin{corollary}
Let $\pi:\mathscr{X}\to\Delta$ be a degeneration of smooth projective
manifolds $\mathscr{X}_s$, $s\in \Delta^\times$ with geometric genus
$p_g>0$. If $\mathscr{X}_0$ is an irreducible variety with only
Gorenstein canonical singularities then $0\in \Delta$ is at finite distance
to $s\ne 0$ with respect to the quasi-Hodge metric.
\end{corollary}

\begin{proof}
By elementary commutative algebra $\mathscr{X}_0$ is Gorenstein implies that
$\mathscr{X}$ is Gorenstein. Moreover $K_{\mathscr{X}_s}=
(K_\mathscr{X}+\mathscr{X}_s)|_{\mathscr{X}_s}=K_\mathscr{X}|_{\mathscr{X}_s}$
for all $s\in\Delta$, hence the theorem on semi-continuity implies that
$h^0(\mathscr{X}_0,K_{\mathscr{X}_0})\ge p_g$.\par

Now for $\mathscr{X}'\to \mathscr{X}$ a resolution of singularities such that
$\mathscr{X}'_0=\bigcup_{i=0}^N X'_i$ is a (not necessarily reduced)
normal crossing divisor, there exists
a component $X'_0$ such that $\phi:X'_0\to \mathscr{X}_0$ is a resolution of
singularities. By definition of canonical singularities we thus have
$K_{X'_0}=\phi^*K_{\mathscr{X}_0} + E$ for $E$ an effective exceptional
divisor. By pulling back canonical sections via $\phi^*$ we conclude that
$h^0(X'_0, K_{X'_0})\ge h^0(\mathscr{X}_0,K_{\mathscr{X}_0})\ge p_g$.\par

Clearly this inequality holds true for $\mathscr{X}'\to\Delta$ a semi-stable reduction
of $\mathscr{X}\to\Delta$. So by the Invariant Cycle Theorem we must have
that $p_g(\mathscr{X}_s)=\sum_{i=0}^N p_g(X_i)$. The above theorem then
implies that $0\in\Delta$ is at finite quasi-Hodge distance.
\end{proof}

\begin{remark}
A normal $\mathbb{Q}$-Gorenstein variety $X$ is said to have (at most)
canonical singularities if for a (hence for any) resolution of singularities
$\phi:Y\to X$, $K_Y=_{\mathbb{Q}}K_X +\sum_i e_i E_i$ with $e_i\ge 0$, where
the sum is over all exceptional divisors. In dimension two, canonical
singularities are precisely $\mathbb{C}^2/G$ for a finite sub-group $G\subset
{\rm SL}(2,\mathbb{C})$, the so called ADE singularities. In dimension three
there is a classification theory due to Reid and Mori \cite{reid}. The generic
hyperplane section of a canonical singularities is again canonical, so the
generic surface slices of canonical singularities are ADE singularities.
\end{remark}

\section{The $m$-th Quasi-Hodge Metrics $g_m$}

Following Viehweg \cite{viehweg}, we will assume that $\mathscr{X}_s$ has
semi-ample canonical bundle for the smooth family $\pi:\mathscr{X}\to S$. This
implies that $\pi_*K_{\mathscr{X}/S}^m$ is a locally free sheaf for each
$m\in\mathbb{N}$ (see Theorem 8.16 in \cite{viehweg}, for this one needs only
that $\mathscr{X}_s$ has Gorenstein canonical singularities). Unlike the case
$m=1$, we do not have an immediate natural hermitian metric on the underlying
vector bundle. Let $m\in\mathbb{N}$ be sufficiently large so that there exists
a divisor $\mathscr{D}\in |K_{\mathscr{X}/S}^m|$ smooth over $S$. By taking an
$m$-cyclic cover along $\mathscr{D}$ (c.f.\ \cite{viehweg} Lemma 2.3), we get a
family
$$ \xymatrix{ \mathscr{Y}\ar[rr]^\phi\ar[rd]_{\tau}&
&\mathscr{X}\ar[ld]^{\pi}\\ & S}
$$
with eigenspace decomposition
$\phi_*\mathscr{O}_{\mathscr{Y}}=\bigoplus_{k=0}^{m-1}(K^{-k}_{\mathscr{X}/S})$
and $R^i\phi_*=0$ for $i>0$ (since $\phi$ is finite). Also from the generalized
Hurwitz formula $K_{\mathscr{Y}/S}=\phi^*K_{\mathscr{X}/S}+(m-1)\mathscr{H}$
with $m\mathscr{H}=\phi^*\mathscr{D}=\phi^*(m K_{\mathscr{X}/S})$, we see that
$\phi^*K_{\mathscr{X}/S}^m=K_{\mathscr{Y}/S}$. Simple spectral sequence
argument and projection formula then show that
$$
\tau_*K_{\mathscr{Y}/S}=\pi_*\phi_*K_{\mathscr{Y}/S}
=\pi_*\phi_*(\phi^*K_{\mathscr{X}/S}^m)
=\pi_*\left(\bigoplus_{k=0}^{m-1}K_{\mathscr{X}/S}^{m-k}\right)
=\bigoplus_{k=1}^{m}\pi_*K_{\mathscr{X}/S}^k.
$$

\begin{defn}
The $m$-th quasi-Hodge metric $g_m$ on $S$ is defined to be
the quasi-Hodge metric attached to $\tau_*K_{\mathscr{Y}/S}$ as defined in \S2.
\end{defn}

Let $S=\Delta^\times$. By Theorem 2.1 (or 2.3), the K\"ahler metric $g_m$ is
incomplete at $0$ if and only if all $H^{n,0}(\mathscr{Y}_s)$ ($s\ne 0$) are
$T^r$-invariant. In view of the above splitting, this should indicate certain
extension properties of pluri-canonical forms in
$H^0(\mathscr{X}_s,K^k_{\mathscr{X}_s})$ ($s\ne 0$) to the central fiber. Here
is our proposal:
\begin{itemize}
\item[I.] Take a semi-stable model $\mathscr{X}\to \Delta$ with $\mathscr{X}_0
=\bigcup_{i=0}^N X_i$. The incompleteness of $g_m$ should be equivalent to
certain twisted pluri-genera equalities
$P_k(\mathscr{X}_s)=\sum_{i=0}^N \tilde{P}_k(X_i)$ for $1\le k\le m$
(c.f.\ Lemma 4.1).\par
\item[II.] The twisted pluri-genera equalities for all $m\in\mathbb{N}$ should
force that there is only one component, say $X_0$, to have non-zero
pluri-genera. Mori theory implies that all the other components are uni-ruled.
It should be enough to verify the equality up to a sufficiently large $m$ which
depends only on $\dim \mathscr{X}_s$. \par
\item[III.] Finally we expect
to contract all $X_i$ with $i\ne 0$ using Mori's extremal contractions. If we
apply directly the MMC, we should then arrive at a model $\mathscr{X}'\to
\Delta$ such that $\mathscr{X}'_0$ has only canonical singularities.
\end{itemize}

To make sense of this, we first show that the incompleteness of $g_m$ is a
necessary condition for degenerations with canonical singularities.
\par

\begin{proof} (of  Proposition 1.3.)
The only key point is that
the $m$-th cyclic covering construction along $D\in|K_X^m|$ can be carried
over to families $\pi:\mathscr{X}\to\Delta$ when $\mathscr{X}_0$ has
canonical singularities of index dividing $m$,
and the resulting covering families $\mathscr{Y}\to\Delta$ has
the property that $\mathscr{Y}_0$ has only Gorenstein canonical singularities.
Hence by Theorem 1.1 (or rather Corollary 2.4) that $g_m$ is incomplete.
\end{proof}

\begin{remark}
Indeed one has the constancy of $P_m$ in $s\in \Delta$ for canonical degenerations
(with semi-ample $K$) and hence the $P_m$ equalities for any semi-stable model. It
is thus natural to conjecture that the $P_m$ equalities for a finite number of $m$'s
are equivalent to the existence of canonical degenerations. Notice that the $P_k$
equality for all $k\le m$ is a stronger assumption than the incompleteness of $g_m$.
\par

It is long conjectured that for smooth projective families $\pi:\mathscr{X}\to
S$ the pluri-genera is locally constant in $s\in S$. Recently Siu \cite{siu}
proved this for smooth varieties of general type. This has then been extended
to families with canonical singularities \cite{kawamata1}. We have chosen to
work on families with semi-ample $K$ to avoid the technicality involved in
Siu's theorem.\par
\end{remark}

In the following two sections, we justify our proposal by verifying it
for curves and by showing that MMC implies the case of Calabi-Yau manifolds.

\section{The Case of Curves}

We will consider curves of genus $\ge 2$. The case of elliptic curves is
easier, which is also a special case of Calabi-Yau manifolds to be discussed
later.

The natural algebro-geometric way to treat this problem is to start with a
semistable degeneration $\mathscr{X}\to \Delta$ with $\mathscr{X}_0=
\bigcup_{i\in I} X_i$ and then try to simplify it through birational
modifications. Since the total space $\mathscr{X}$ is a surface, the
modifications needed are simply contraction of $(-1)$ curves. So we may assume
that $\mathscr{X}\to S$ is relative minimal and semi-stable and
$g(\mathscr{X}_s)\ge 2$. Let $g_i = g(X_i)$ and $d_i=\sum_{j\ne i}X_j.X_i =
(\sum X_j - X_i).X_i = -X_i^2$. The famous stable reduction theorem for curves
states that we may further contract $(-2)$ curves so that every component $X_i$
is stable in the sense that if $g_i=0$ then $d_i\ge 3$ (and if $g_i=1$ then
$d_i\ge 1$, which is always true here since $\mathscr{X}_0$ is connected). The
only subtle point is that $\mathscr{X}$ may have $A_n$ type singularities.
Since this will not affect our later discussion, for the sake of simplicity we
will assume we are already in a stable reduction.

For $m\in \mathbb{N}$ such that $K_{\mathscr{X}/S}^m$ is $S$-very ample, let
$\mathscr{D}\in|K^m_{\mathscr{X}/S}|$ be a smooth member which does not contain
any special point in any $X_i\cap X_j$ and the singular points of
$\mathscr{X}$. We also assume that $S$ is small enough so that $S\cong \Delta$
and $K_S\cong \mathscr{O}_S$. By Theorem 1.1, the quasi-Hodge metric $g_m$
constructed from the $m$-th cyclic cover $\mathscr{Y}\to\mathscr{X}$ along
$\mathscr{D}$ is incomplete if and only if $p_g(\mathscr{Y}_s)=\sum_{i\in
I}p_g(Y_i)$. This is because $\mathscr{Y}\to S$ is a stable degeneration by
construction and the presence of $A_n$ singularities does not affect the result
--- further blowing-ups gives $(-2)$ curves which have no contribution to the
geometric genus. By the generalized Hurwitz formula, we have seen in \S3 that
$p_g(\mathscr{Y}_s)=\sum_{k=1}^m P_k(\mathscr{X}_s)$. Moreover, we have
\begin{lemma}
$p_g(Y_i)=\sum_{k=1}^m h^0(X_i, \tilde{K}_{X_i}^k)$ where $\tilde{K}_{X_i}^k$ is the
twisted (or logarithmic) pluri-canonical sheaf defined as
$K_{X_i}^k\otimes\mathscr{O}_{X_i}(\sum_{j\ne i} X_j\cap X_i)^{k-1}$.
\end{lemma}

\begin{proof}
Let $\tau:\mathscr{Y}\to S$ be the new families with
$\mathscr{Y}_0=\bigcup_{i\in I}Y_i$. By the construction, $\phi_i:Y_i\to X_i$
is an $m$-cyclic cover along
$$
\mathscr{D}|_{X_i}=mK_\mathscr{X}|_{X_i}=m(K_{X_i}-X_i|_{X_i})
=m(K_{X_i}+\sum\nolimits_{j\ne i}X_j\cap X_i)=:mD_i.
$$
Then we have eigenspace decomposition ${\phi_i}_*\mathscr{O}_{Y_i}
=\bigoplus_{k=0}^{m-1}\mathscr{O}(D_i)^{-k}$. The same
proof as in \S3 gives ${\phi_i}_* K_{Y_i}=K_{X_i}\otimes
\bigoplus_{k=0}^{m-1}\mathscr{O}(D_i)^{k}$.
\end{proof}

\begin{proof}(of Theorem 1.4.)
We will show that if $g_m$ is incomplete for three
$m$'s then in a stable reduction there is only one
component $X_0$ in $\mathscr{X}_0$,
which is then smooth of genus $g(X_0)=g(\mathscr{X}_s)$.

From the Riemann-Roch formula,
\begin{eqnarray*}
h^0(\tilde{K}_{X_i}^k)-h^1(\tilde{K}_{X_i}^k)
  \!\!\!&=&\!\!\! k(2g_i-2)+ (k-1)d_i + (1-g_i)\\
         &=&\!\!\! (2k-1)(g_i -1) + (k-1)d_i
\end{eqnarray*}

It is clear that $h^1(\tilde{K}_{X_i}^k)=0$ for $k\ge 2$ by the stability
condition. Also $h^1(\tilde{K}_{X_i})=h^0(\mathscr{O}_{X_i})=1$. So the
equality
$$
\sum\nolimits_{k=1}^m h^0(K_{\mathscr{X}_s}^k)=\sum\nolimits_{i\in
I} \sum\nolimits_{k=1}^m h^0(\tilde{K}_{X_i}^k)
$$
becomes
$$
1 + \sum\nolimits_{k=1}^m (2k-1)(g-1)=|I| + \sum\nolimits_{i\in I}
\sum\nolimits_{k=1}^m \big[(2k-1)(g_i -1) + (k-1)d_i\big]
$$
which is
$$
m^2(g-\sum\nolimits_{i\in I}g_i)= -(|I|-1)m^2 + \frac{m(m-1)}{2}
\sum\nolimits_{i\in I}d_i + (|I|-1).
$$
If this is true for any three values of $m$'s then we get $|I|=1$, say
$I=\{0\}$ and $g=\sum_{i\in I} g_i = g(X_0)$ as desired.
\end{proof}

\section{The Case of Calabi-Yau Manifolds}

We start with the following corollary of Theorem 1.1.

\begin{corollary}
Let $\mathscr{X}\to \Delta$ be a semi-stable degeneration of Calabi-Yau
manifolds. Then $g_H$ is incomplete at $0\in\Delta$
if and only if there is an irreducible component
$X_i \subset \mathscr{X}_0$ such that $h^{n,0}\ne 0$. This is equivalent to that
there is exact one component with $h^{n,0}=1$.
\end{corollary}

\begin{lemma}
For any relative minimal model $\mathscr{X}\to \Delta$ of a degeneration of
Calabi-Yau manifolds $\mathscr{X}'\to \Delta$, if $\mathscr{X}_0 =\sum_{i=0}^N X_i$
has more than one component then each $X_i$ has $-K_{X_{i,{\rm red}}}$ a nontrivial
effective divisor on $X_{i,{\rm red}}$.
\end{lemma}

\begin{proof}
For divisors we use the notation ``$\sim$'' to denote $\mathbb{Q}$-linear
equivalence. Since $\pi:\mathscr{X}\to \Delta$ is a holomorphic function, we have
that $\mathscr{X}_s\sim 0$ for any $s\in\Delta$. Also if $D\subset \mathscr{X}$ is a
divisor such that $\pi(D)$ is not a point, then $D|_{\mathscr{X}_s}$ will be a
nontrivial divisor of $\mathscr{X}_s$.\par

Since $\mathscr{X}_s$ is birational to $\mathscr{X}'_s$, it is a terminal Calabi-Yau
variety for $s\ne 0$, that is $K_{\mathscr{X}_s}=0$. By adjunction formula on such
$\mathscr{X}_s$, we conclude that $K_{\mathscr{X}}$ is supported on the central
fiber and so is of the form $K_{\mathscr{X}}=\sum a_i X_i$ with $a_i\in\mathbb{Q}$.
\par

Since $\sum X_i =\mathscr{X}_0 \sim 0$, we may adjust $a_i$ so that $\max\, a_i =
0$. Let $I=\{i\,|\,a_i = 0\}$. For $i\in I$ and a curve $\ell\subset X_i$ one has $$
K_{\mathscr{X}}.\ell = \sum\nolimits_{j\not\in I} a_j(X_j.\ell)\le 0. $$ If $I\ne
\{0,\ldots,N\}$, we may choose $\ell$ such that $X_j.\ell > 0$ and $a_j < 0$ by the
connectedness of $\bigcup X_i$. But then $K_{\mathscr{X}}.\ell < 0$, contradicts to
the nefness of $K_{\mathscr{X}}$.\par

So we must have $a_i = 0$ for all $i$. That is, $K_{\mathscr{X}}=0$. Now take any
component $X_i$, one has
\[
K_{X_i}=(K_{\mathscr{X}}+X_i)\big|_{X_i} = X_i|_{X_i}=-\sum\nolimits_{j\ne i}
X_j|_{X_i},
\]
which by the connectedness again is a nontrivial negative effective
divisor if there are more than one components. The case for $X_{i,{\rm red}}$ is
entirely similar.
\end{proof}

\begin{proof} (of Proposition 1.2.) We first take a semi-stable reduction
$\mathscr{X}'\to \Delta$ of our original degenerations. This may require a finite
base change on the base, but this will not change the finite distance condition on
the metric. We want to conclude that any relative minimal model (if there exists
any) of it, denoted by $\mathscr{X}\to \Delta$ has only one component in the central
fiber. Let $\mathscr{X}'_0=\bigcup_{i=0}^N X'_i$ with $X'_0$ the unique component
with a canonical section $\Omega\in\Gamma(X'_0,K_{X'_0})$.\par

By the conjectural construction of relative minimal models, one uses only divisorial
contractions and flips of relative extremal rays, hence only the central fiber will
be modified since the general fibers are already smooth Calabi-Yau. Notice that
during the process the proper transform of $X'_0$ is never contracted. Indeed, if a
component $W$ in the central fiber is contracted in some step, then $W$ is covered
by (extremal) rational curves \cite{kawamata2}. This will continue to hold true for
any smooth model $Y$ of $W$ hence $\kappa(Y)=-\infty$ and so $Y$ can not be
$X'_0$.\par

By Lemma 5.2, if there are more than one components in $\mathscr{X}_0$ then any
component $X_i$ has $K_{X_i}=-D\ne 0$. If $X_i$ is not normal, passing to
normalization $\psi:W\to X_i$ with conductor $C\subset W$ can only make
$K_W=\psi^*K_{X_i}-C$ more non-effective. Hence passing to any smooth model
$\phi=\phi'\circ\psi:Y\to X_i$ with $\phi':Y\to W$, $K_Y=\phi'^*K_W + E$ shows that
$K_Y=-\phi^* D -\phi'^*C + E$. This is never an effective divisor because $E$ is
$\phi'$-exceptional. By the birational invariance of $p_g$ among smooth models, this
contradicts to the existence of $\Omega$ on $X'_0$. Hence $\mathscr{X}_0$ has only
one component and
$K_{\mathscr{X}_0}=(K_{\mathscr{X}}+\mathscr{X}_0)|_{\mathscr{X}_0}
=\mathscr{X}_0|_{\mathscr{X}_0}=0$.\par

To see that $X:=\mathscr{X}_0$ has at most canonical singularities, let $\phi:Y\to
X$ be a resolution of singularities. Clearly $\Gamma(Y,K_Y)\cong\mathbb{C}$ since
$X$ has a smooth model with $p_g=1$. In particular, $K_Y$ is effective. Also
$\Gamma(X,K_X)\cong\mathbb{C}$ since $K_X=0$. If $X$ is normal, since it is
Gorenstein we have $K_Y=\phi^*K_X + E = E$ which shows that $E$ is effective and so
$X$ has at most canonical singularities. If $X$ is not normal, let $\psi:W\to X$ be
the normalization with non-zero conductor divisor $C\subset W$ and let $\phi$
factors through $\psi$ as $\phi=\phi'\circ\psi$ with $\phi':Y\to W$. Then
$K_W=\psi^*K_X - C = -C$ and $K_Y = \phi'^*K_W + E'=-\phi'^*C +E'$, which is never
effective, a contradiction. Hence $X$ is normal and is in fact a Calabi-Yau variety
with at most canonical singularities. The proof is thus completed.
\end{proof}

\section{Essential Incomplete Boundary Points}

A degeneration over the unit disk is called {\it non-essential} if it admits a
finite base change so that the punctured family can be completed into a smooth
family. Otherwise it is called {\it essential}. Clearly a necessary condition of a
degeneration to be non-essential is that the monodromy is of finite order, or $N=0$.
In particular it is of finite distance. In fact for K3 surfaces any finite distance
degeneration is non-essential \cite{kulikov}. It is thus interesting to see whether
there exists essential finite distance degenerations in higher dimensions.\par

The classical Picard-Lefschetz theory states that for a Lefschetz pencil (i.e.\
nodal degenerations) each node, or ODP, $p_i$ will correspond to a vanishing cycle
$\sigma_i\in H_n(\mathscr{X}_s,\mathbb{Z})$ and all these $\sigma_i$'s generate the
space of vanishing cycles $V$ which is the kernel of the map $H_n(\mathscr{X}_s)\to
H_n(\mathscr{X}_0)=H_n(\mathscr{X})$. Moreover the monodromy $N\ne 0$ if and only if
$V\ne 0$. However, there is no general local criterion to test whether $\sigma_i\ne
0$ for a given ODP $p_i$. The $\sigma_i$ is always trivial for $n$ even and must be
trivial even for $n$ odd if $H_n(\mathscr{X}_s)=0$ --- for example, nodal
degenerations for cubic or quartic threefolds.\par

Here we show that $V\ne 0$ for any nodal degenerations of Calabi-Yau 3-folds.
\footnote{The problem on essential finite distance boundaries is discussed in length
in \cite{wang}, \S3-\S4. However, the argument in \cite{wang}, \S3 concerning nodal
degenerations is not complete. The author is grateful to M.\ Gross for discussions
on this issue.}

\begin{proof} (of Theorem 1.5, part (1))
First of all, a nodal threefold $\mathscr{X}_0$ always admits (not necessarily
projective) small resolutions $X\to \mathscr{X}_0$ with smooth rational curves
$X\supset C_i\to p_i\in \mathscr{X}_0$ contracted to ODP's. In the case of
Calabi-Yau threefolds (Gorenstein threefolds with trivial canonical bundle) with
$h^1(\mathscr{O})=0$, the existence of global smoothing $\mathscr{X}\to\Delta$ of
$\mathscr{X}_0$ forces that there are nontrivial relations of $[C_i]\in H_2(X)$ by
Friedman's result \cite{freidman1}, \cite{freidman2}. That is, the canonical map
$e:\bigoplus_i \mathbb{Z}[C_i]\to H_2(X,\mathbb{Z})$ has nontrivial kernel dimension
$\rho >0$. Consider the resulting surgery diagram:
\begin{eqnarray*}
&X&\\
&\downarrow&\\
&{}\,\mathscr{X}_0\!&\!\!\!\subset\mathscr{X}\supset \mathscr{X}_s
\end{eqnarray*}
It has the following local description: let $V_i\ni p_i$ be a contractible
neighborhood of an ODP, $V'_i\subset \mathscr{X}_s$ be the smoothing of $V_i$
and $U_i\subset X$ be the inverse image of $V_i$. Then
\begin{itemize}
\item[I.] $U_i$ is a deformation retract neighborhood of $C_i$ and so has the
homotopy type of $S^2\sim D^4\times S^2$.
\item[II.] $V'_i$ has the homotopy type of $S^3\times D^3$. Where
the sections $\sigma_i\sim S^3$ are the vanishing cycles.
\item[III.] The surgery from $X$ to $\mathscr{X}_s$ is induced from
$\partial(D^4\times S^2)=S^3\times S^2=\partial(S^3\times D^3)$.
\end{itemize}

Let us assume that there are $k$ ODP's. An immediate consequence
is the Euler number formula:
$$
\chi(X)-k\chi(\mathbb{P}^1)=\chi(\mathscr{X}_0)-k\chi({\rm pt})=
\chi(\mathscr{X}_s)-k\chi(S^3).
$$

Let $W$ be the ``common open set'' of $X$, $\mathscr{X}_0$ and $\mathscr{X}_s$ away
from all points $p_i$'s such that $W$ and $V_i$'s cover $\mathscr{X}_s$ etc.
A portion of the Mayer-Vietoris sequence of the covering
$\{W,\,V'_i\}$ of $\mathscr{X}_s$ gives
$$
0\to H_3(W)\to H_3(\mathscr{X}_s)\to \bigoplus\nolimits_i \mathbb{Z}[C_i]\to
H_2(X)\to H_2(\mathscr{X}_s)\to 0.
$$
Hence that $b_2(X)=b_2(\mathscr{X}_s) + (k - \rho)$.\par

Take into account of $b_2(\mathscr{X}_0)=b_2(\mathscr{X}_s)$ and
$b_4(\mathscr{X}_0)=b_4(X)$ (which also follows from suitable Mayer-Vietoris
sequences), simple manipulations show that
$b_3(\mathscr{X}_s)=b_3(\mathscr{X}_0) + \rho$. Comparing with the
(Mayer-Vietoris) sequence defining the vanishing cycles: $$
\bigoplus\nolimits_i \mathbb{Z}[\sigma_i]\to H_3(\mathscr{X}_s)\to
H_3(\mathscr{X}_0)\to 0, $$ we conclude that $\rho$, which is non-zero,
is the dimension of $V$.
\end{proof}

It is possible for
an essential degeneration to have trivial monodromy. In \cite{freidman2},
Friedman remarked that Clemens has constructed families of quintic
hypersurfaces in $\mathbb{P}^4$ acquiring an $A_2$ singularity and have
monodromy of finite order. He then asked whether this family can be filled in
smoothly up to a base change. This has been
answered negatively in \cite{wang}. We recall the statement here for the reader's
convenience.
\begin{theorem}
Let $\mathscr{X}\to\Delta$ be a projective smoothing
of a nontrivial Gorenstein terminal minimal threefold $\mathscr{X}_0$ over the unit disk. Then,
up to any finite base change,
$\mathscr{X}\to\Delta$ is not $\Delta$-birational to a projective smooth family
$\mathscr{X}'\to\Delta$ of minimal threefolds.
\end{theorem}
In view of Proposition 1.3, this implies that terminal
degenerations of smooth minimal threefolds
provide essential incomplete boundary points of the moduli spaces
with respect to the quasi-Hodge metrics $g_m$ for all $m$. This proves
Theorem 1.5, part (2).

\section{Remarks on Moduli and Weil-Petersson Metrics}

Canonical singularity naturally occurs in minimal and canonical models in algebraic
geometry. It also plays a significant role in string theory through the connection
with Calabi-Yau manifolds. On the other hand, in various situations it behaves just
like smooth points. A nice example is Viehweg's program on constructing
quasi-projective moduli spaces of polarized manifolds \cite{viehweg}. He showed that
it is essentially the same proof to include varieties with canonical singularities
as long as the deformation invariance can be verified. This was recently proved by
Kawamata \cite{kawamata1} based on \cite{siu}.\par

Along different lines, the author had attempted to understand the boundary of moduli
spaces of Calabi-Yau manifolds from the differential geometric point of view
\cite{wang}. It was found that the natural Weil-Petersson metric on the moduli space
is incomplete, therefore the metric completion of moduli spaces becomes an important
problem. It was proved that degenerations of Calabi-Yau manifolds with at most
canonical singularities are at finite Weil-Petersson distance. It was also
conjectured there that the converse holds. Its truth would imply that Viehweg's
enlarged moduli spaces coincide with the metric completion of the moduli spaces,
hence a perfect match between viewpoints in algebraic geometry and differential
geometry. Now this follows from the minimal model conjectures by Proposition
1.2.\par

The Weil-Petersson metric on the moduli space of Calabi-Yau manifolds is defined as
the variation of the underlying Ricci flat metrics. For a given polarized Calabi-Yau
family $\mathscr{X}\to S$ with Ricci flat metrics $g(s)$ on $\mathscr{X}_s$, under
the Kodaira-Spencer map $\rho: T_{S,s}\to H^1(\mathscr{X}_s,T_{\mathscr{X}_s})\cong
\mathbb{H}^{0,1}_{\bar\partial}(T_{\mathscr{X}_s})$ (harmonic forms with respect to
$g(s)$), we have for $v$, $w\in T_s(S)$,
$$
g_{WP}(v,w):=\int_{\mathscr{X}_s}\langle\rho(v),\rho(w)\rangle_{g(s)}.
$$
It is a
quite surprising fact that $g_{WP}$ admits a Hodge theoretic description. Indeed
$g_{WP}=g_H$. This follows from the fact that the holomorphic volume form
$\Omega(s)$ is parallel with respect to $g(s)$, which again is equivalent to the
Ricci flat condition
$$
\Omega(s)\wedge\bar\Omega(s) = f(s)\,\omega_{g(s)}^n
$$
for
a constant $f(s)$ depending only on $s\in S$. Indeed $f(s)$ is the point-wise length
square of $\Omega(s)$ if we normalize the volume to be 1. This viewpoint provides an
alternative differential geometric way to look at the above canonical singularity
conjecture without using the minimal model theory.

With the $\Omega$ chosen as in \S2, the incompleteness of $g_H=g_{WP}$ of a
punctured Calabi-Yau family $\pi:\mathscr{X}\to \Delta^\times$ is equivalent to the
continuity of $f(s)$ over $0\in \Delta$. We attempt to show from this the uniform
boundedness of diameter of $\mathscr{X}_s$ for all $s\in \Delta^\times$. With this
done, we may then proceed by using the theory of Hausdorff convergence. The details
of this differential geometric approach will appear in a separate work.
\par

\bibliographystyle{amsplain}

\end{document}